\documentclass[11pt]{amsart}
\usepackage{geometry}                
\usepackage{graphicx}
\usepackage{amssymb}
\usepackage{epstopdf}
\title{Group Theory and Modern Dance Composition}
\author{Ashley Walls White}
\date{May 23, 2020}                                           

\begin{document}
\maketitle

\section*{Abstract}
This paper will examine the spatial reference systems typically used in Laban Movement Analysis (LMA), and the consequences of group actions on these systems. The elementary notions of inversion and transposition in choreographic composition can be defined in such a way that they can be shown to be group homomorphisms in all the reference systems of LMA. The notions of orbits and stabilizers on polyhedra are used to mathematically define these choreographic devices, and these same notions can be used to define new choreographic devices on the standard polyhedra used for spatial reference in dance.

\section*{Introduction}
Laban Movement Analysis (LMA) was developed by Rudolf Laban as a means of describing and documenting human movement, and is used as a tool for studying and analyzing modern dance. The main aspects of LMA I will examine are the use of regular polyhedra as spatial reference systems, and the development of movement scales within these figure. The cube, octahedron, and icosahedron are used to determine 26 distinct spatial directions, each coming from the vertices of these polyhedra. The notions of group actions, or symmetries, on these figures are surprisingly useful in modern dance choreography. In choreographic practice, as in music composition, a thematic phrase of movement material is developed, and then choreographic devices are applied to this phrase to generate related and meaningful iterations of new movement material. One particular devices is the notion of inversion. Informally, inversion is defined as taking a movement to its spatial opposite, for instance inversion could take what was in the front space of the body to the backspace, or an action with an upright directional intent to having a downward directional intent. Inversion can be formally defined in the spatial reference systems, and we can show that this notion can be defined in terms of a group action, if we take the group to be the vertices of any of the three noted polyhedrons.

\section*{Inversion}

\subsection*{The Octahedron}
We will start by examining the octahedron. We can consider the body situated inside an octahedron. Let $v_1$ be the head, $v_2$ be the feet, and $v_4$ be the front of the body so that $v_6$ is the middle side left direction and $v_5$ is the middle side right. We define inversion as taking the directional intent of an action to its spatial opposite. Using the vertices of the octahedron, we would hope that the inversion of an action in the direction of $v_4$ would be an action in the direction of $v_3$, and similarly the inversion of an object is the direction of $v_1$ would be the same action in the direction $v_2$.\\
\\

\begin{figure}[h!]
  \centering
  \includegraphics[scale=0.35]{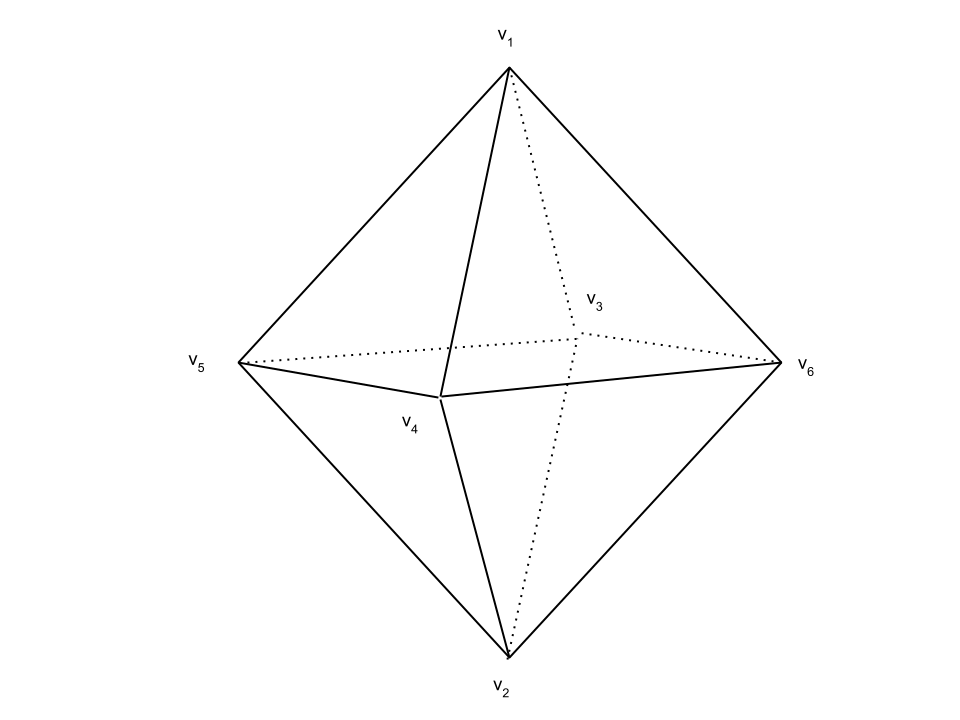}
  \caption{Octahedron}
\end{figure}

Formally consider the direction of an action to be the set containing the vertex that defines a particular direction, so an action that takes place in the middle right direction would be in the direction $\{v_6\}$. \\
\\

Consider the set of the vertices $V = \{v_1, v_2, v_3, v_4, v_5, v_6\}$ and the group of rotational symmetries on the octahedron. Considering $\{v_i\}$ for $i = 5, 6$, the stabilizer of $v_i$ is the set of rotations that leave $v_i$ fixed. The stabilizer $H = stab(v_i)$ fixes both $v_5$ and $v_6$ and induces a group action on the octahedron. There is only one permutation other than the identity in the stabilizer of $v_i$. \\

\indent$H = stab(v_i) = \{e, (v_1, v_2)(v_3, v_4)\}$, for $i = 5,6$\\
\\
We can consider the orbits under this group action. \\
\\
\indent $orb(v_1) = \{v_1, v_2\}$\\
\indent $orb(v_3) = \{v_3, v_4\}$  \\
\\
Since the orbits partition the set $V$, for an action in any direction, we can define the inversion of that action to be the only other element in the action's orbit. Considering the action in the middle back direction, $\{v_3\}$, we can define the inversion\\
\\
\indent$I(v_3) = orb(v_3) - \{v_3\} = \{v_4\}$\\
\\
which gives us an action in $\{v_4\}$, the middle front direction, which is what we expected inversion to be based on our informal definition. Similarly, \\
\\
\indent$I(v_4) = orb(v_4) - \{v_4\} = \{v_3\}$\\
\indent$I(v_1) = orb(v_1) - \{v_1\} = \{v_2\}$\\
\indent$I(v_2) = orb(v_2) - \{v_2\} = \{v_1\}$\\
\\
Actions in the mid level directly left and directly right are inversions of each other, and taking the stabilizer of any other vertex in the octahedron puts these directions in the same orbit.\\
\\
In general, for an invertible action $\{v_i\}$, define \\
\\
\indent$I(v_i) = orb(v_i) - \{v_i\} $\\
\\
Inversion in the octahedron sets the stage for defining choreographic devices on the icosahedron, which will provide a richer description of directions in space. In the icosahedron we will continue the practice of examining orbits and stabilizers to define inversion, but a clearer distinction will be made between group actions that define each type of inversion, high-low, front-back, and left right.  
\subsection*{The Icosahedron}For our purposes it will be most convenient to consider the vertices of the icosahedron formed by the corners of the mutually orthogonal intersecting rectangles. In anatomy we consider the vertical, horizontal, and sagittal planes on the body. The vertical plane cuts the body into front and back halves, the horizontal plane cuts the body into top and bottom halves, and the sagittal plane cuts the body into left and right halves. \cite{bart}  \\

\begin{figure}[h!]
  \centering
  \includegraphics[scale=0.35]{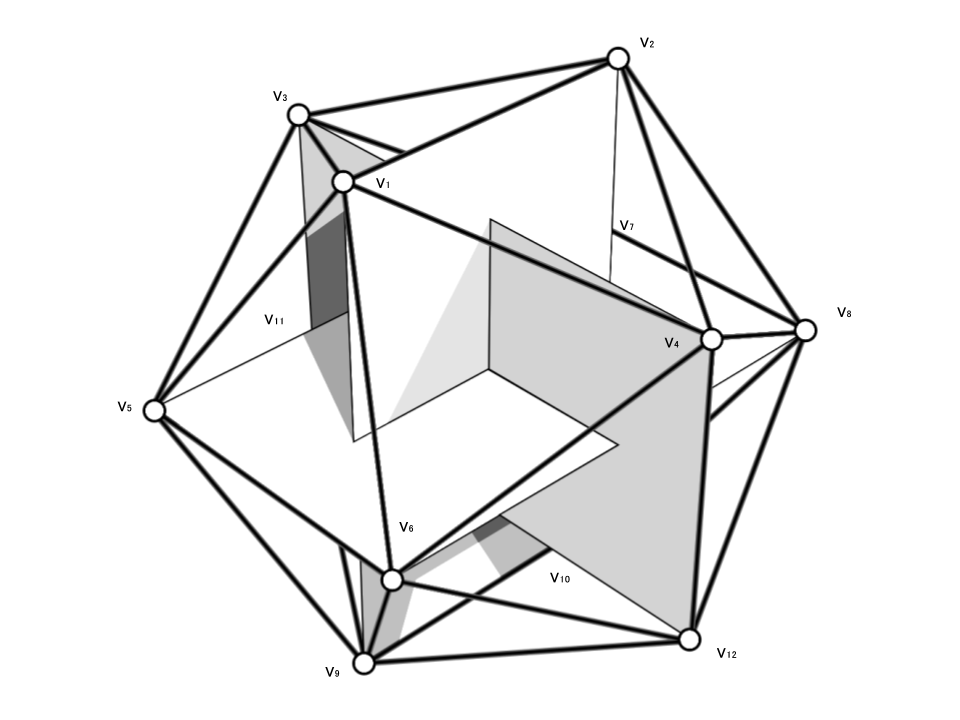}
  \caption{Icosahedron from Orthogonal Rectangles}
\end{figure}

In the figure above, the blue plane is the sagittal plane, the green plane is the horizontal plane, and the red plane is the vertical plane. If we situate the body inside the icosahedron so that the head is at the top of the sagittal plane, the vertices correlate to the following directions of the body in space.\\
\\
\indent \begin{tabular}{ll}
$v_1$ - forward high & $v_7$ - middle right back \\
$v_2$ - back high & $v_8$ - middle left back\\
$v_3$ - high right & $v_9$ - forward low \\
$v_4$ - high left & $v_{10}$ - back low \\
$v_5$ - middle right forward & $v_{11}$ - low right \\
$v_6$ - middle left forward & $v_{12}$ - low left\\
\end{tabular}\\
\\
\\
\\
Let's first begin to consider the group of rotational symmetries on the figure, $I_r$. We only want symmetries that will keep the same orientation of the head and feet in the sagittal plane. Allowing $I_r$ to act on the edges of the icosahedron, note that keeping the "head" in place might correlate to taking the stabilizer of the top edge connecting $v_1$ and $v_2$, call it $e_1$. Considering $stab(e1)$ we find that this is just the rotational symmetry $\rho_{180}$.\\
\\
The rotational symmetry $\rho_{180}$ is working about a 2-fold axis that passes through the midpoint of $e_1$ to the midpoint of the opposite edge, the edge connecting $v_9$ and $v_{10}$. This symmetry induces a group action on the vertices, from which we can then determine the partition of the vertices into orbits. \\
\\
\\
\indent $orb(v_1) = \{v_1 = e, v_2\}$\\
\indent $orb(v_3) = \{v_3 = e, v_4\}$  \\
\indent $orb(v_5) = \{v_5 = e, v_7\}$\\
\indent $orb(v_6) = \{v_6 = e, v_8\}$  \\
\indent $orb(v_9) = \{v_9 = e, v_{10}\}$\\
\indent $orb(v_{11}) = \{v_{11} = e, v_{12}\}$  \\
\\
From choreographic composition, there are many ways to interpret inversion. The three we will focus on in this section are from forward to back, from low to high and from left to right. For the front-back inversion, the inversion of forward high, $v_1$, is back high, $v_2$, and from above those two vertices are in the same orbit. All the orbits defined above partition vertices into pairs with their front-back inversion, $I_{fb}$, except for the permutations $(v_3, v_4)$ and $(v_{11}, v_{12})$. These two orbits give define a relation between left and right directions, which does not fit into the definition of a front-back inversion, so we need an action on the vertices that will fix this entire plane, i.e. the vertices $v_3, v_4, v_{11}, v_{12}.$\\
\\
Considering the $stab(v_i, v_{i + 1})$, for $i = 3, 11$ we arrive at the following:\\
\\
\indent$H = stab(v_i, v_{i + 1}) = \{e, (v_1, v_9)(v_2, v_{10})(v_3, v_{11})(v_4, v_{12})\}$\\
\\
The non-identity permutation defines the following orbits:\\
\\
\indent $orb(v_1) = \{v_1 = e, v_2\}$\\
\indent $orb(v_5) = \{v_2 = e, v_{7}\}$  \\
\indent $orb(v_6) = \{v_3 = e, v_{8}\}$\\
\indent $orb(v_9) = \{v_4 = e, v_{10}\}$  \\
\\
This partitions the orbits into sets containing an action and that action's front-back inversion. Another symmetry will have to be determined to define low-high inversion and left-right inversion.\\
\\
\subsection*{Low-High Inversion} For low-high inversion, $I_{lh}$, we again cannot employ a rotational symmetry.  First consider again $\rho_{180}$, but consider it about the axis that goes through the horizontal place, entering at the edge connecting $v_5$ and $v_6$, and existing the edge connecting $v_7$ and $v_8$. The partition of the vertices into orbits under this action is as follows:\\
\\
\indent $orb(v_1) = \{v_1 = e, v_9\}$\\
\indent $orb(v_2) = \{v_2 = e, v_{10}\}$  \\
\indent $orb(v_3) = \{v_3 = e, v_{11}\}$\\
\indent $orb(v_4) = \{v_4 = e, v_{12}\}$  \\
\indent $orb(v_5) = \{v_5 = e, v_{6}\}$\\
\indent $orb(v_7) = \{v_{7} = e, v_{8}\}$  \\
\\
This works except for all the vertices in the horizontal plane. These do not have low-high inversions, so we need an action on the vertices that will fix this entire plane, i.e. the vertices $v_5, v_6, v_7, v_8.$\\
\\
Considering the $stab(v_i, v_{i + 1})$, for $i = 5..8$ we get:\\
\\
\indent$H = stab(v_i, v_{i + 1}) = \{e, (v_1, v_9)(v_2, v_{10})(v_3, v_{11})(v_4, v_{12})\}$\\
\\
The non-identity permutation defines the following orbits:\\
\\
\indent $orb(v_1) = \{v_1 = e, v_9\}$\\
\indent $orb(v_2) = \{v_2 = e, v_{10}\}$  \\
\indent $orb(v_3) = \{v_3 = e, v_{11}\}$\\
\indent $orb(v_4) = \{v_4 = e, v_{12}\}$  \\
\\
which correctly determines the sets of elements with their low-high inversion. \subsection*{Left-Right Inversion} Left-Right inversion is similarly defined by taking the stabilizer of sagittal plane. \\
\\
Considering the $stab(v_i, v_{i + 1})$, for $i = 1..2$ we get:\\
\\
\indent$H = stab(v_i, v_{i + 1}) = \{e, (v_3, v_4)(v_5, v_{6})(v_7, v_{8})(v_{11}, v_{12})\}$\\
\\
The non-identity permutation defines the following orbits:\\
\\
\indent $orb(v_3) = \{v_3 = e, v_4\}$\\
\indent $orb(v_5) = \{v_5 = e, v_{6}\}$  \\
\indent $orb(v_7) = \{v_7 = e, v_{8}\}$\\
\indent $orb(v_{11}) = \{v_{11} = e, v_{12}\}$  \\
\\
We find that what is referred to as "inversion" in modern dance choreography follows a reflective symmetry in the icosahedron.
\\
\\
\section{Additional Uses for Group Actions}
\subsection*{Subsequent inversions on the icosahedron}
In choreography an inversion can be applied to an action multiple times in a row. Consider front-back inversion on the icosahedron and actions $A$ and $B$, where action $A = \{v_5\}$ is in the direction of $v_5$ and action $B = \{v_8\}$ is in the direction of $v_8$.\\
\\
\indent $I_{fb}(A) = \{v_7\}$ \\
\indent $I_{fb}(B) = \{v_6\}$\\
\\
Define the operation $P \oplus R$ on actions $P$ and $R$ to be performing an action in the direction of P, and then performing an action in the direction $R$, the resulting action just being in the direction $R$. Then,\\
\\
\indent  $I_{fb}(A) \oplus I_{fb}(B) = \{v_6\}$. \\
\\
And\\
\indent $A \oplus B = \{v_8\}$, \\
\\
so that \\
\indent $I_{fb}(A \oplus B) = \{v_6\}$.\\
\\
Performing two actions in sequence and then taking the inversion yields the same results as doing the inversions of each of the actions in sequence. We can define the inversion of a sequence of events to be sequence of each individual action's inversion. It should be noted that this definition is only concerned with the end destination of the movement, but the not the sequence of movements taken to get there.  \\
\\
\subsection*{Laban's normal zones as group actions}
Laban identified five principle zones of the body associated with the head, each arm, and each leg. \cite{laban} Considering the body positioned in the icosahedron as before, consider the right and left arms at high right $\{v_3\}$ and high left $\{v_4\}$ respectively, with the legs similarly at $\{v_{11}\}$ and $\{v_{12}\}$.  If we take the stabilizer of any of these vertices we can create orbits which roughly correspond to Laban's normal zones for the arms and legs. Consider first $stab(v_3)$. Under this the following orbits are created.\\
\\
\indent $orb(v_1) = \{v_1, v_2, v_8, v_{12}, v_6\}$\\
\indent $orb(v_3) = \{v_3, v_7, v_{10}, v_9, v_5\}$\\  
\indent $orb(v_4) = \{v_4\}$\\
\indent $orb(v_{11}) = \{v_{11}\}$\\
\\
The first, $orb(v_1)$ gives a group of vertices within a normal range for the left arm, with $v_6$ being the standard position for the left arm.  The normal range, which relates to Laban's normal zones, for any arm or leg can be defined as the orbit of any of the vertices adjacent to the vertex that gives the standard position for that limb under the stabilizer of the vertex that gives the standard position. Thus the normal zone for each limb is the following:\\
\\
\hspace*{2cm}
\begin{tabular}{ccc}
Limb & Standard Position & Normal Range \\
\hline
Left Arm & $\{v_4\}$&  $orb(v_1) = \{v_1, v_2, v_8, v_{12}, v_6\}$\\\\
Right Arm & $\{v_3\}$&  $orb(v_1) = \{v_1, v_2, v_7, v_{11}, v_5\}$\\\\
Left Leg & $\{v_{12}\}$&  $orb(v_9) = \{v_6, v_4, v_8, v_{10}, v_9\}$\\\\
Right Leg & $\{v_{11}\}$&  $orb(v_9) = \{v_5, v_3, v_7, v_{10}, v_9\}$\\\\
\end{tabular}

\section*{Laban's Directions in Space}
Laban gave a codified system for representing directions in space. Henceforth we will refer to the directions using his symbols, which are listed below \cite{laban}.
\\
\begin{figure}[!h]
  \centering
  \includegraphics[scale=0.35]{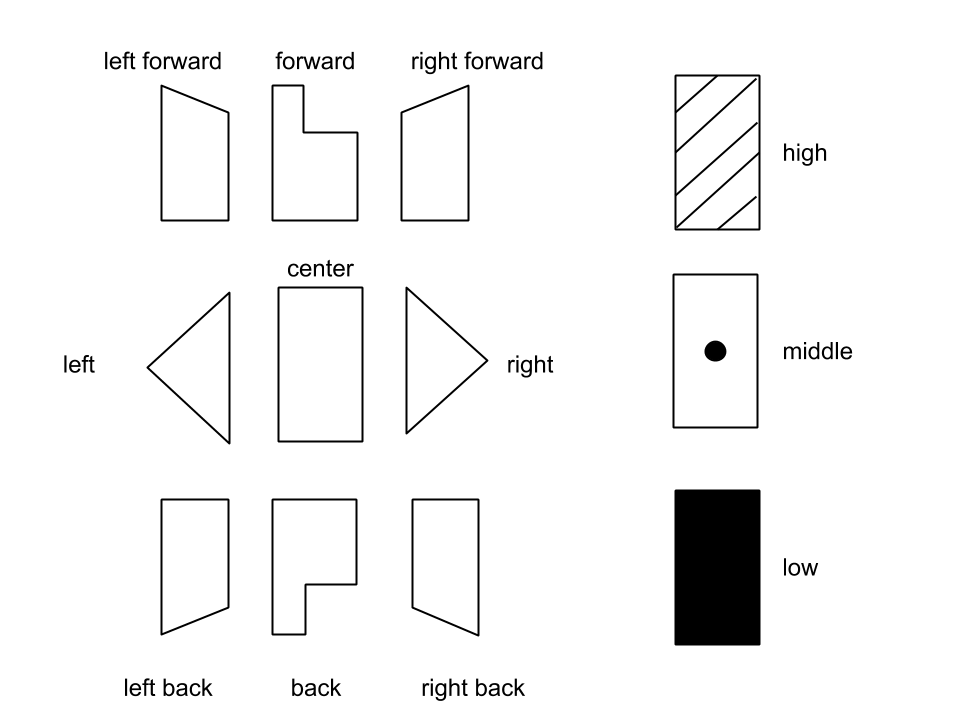}
  \caption{Laban's Directions in Space}
\end{figure}

Below is the icosahedron with these directions labeled.

\begin{figure}[!h]
  \centering
  \includegraphics[scale=0.3]{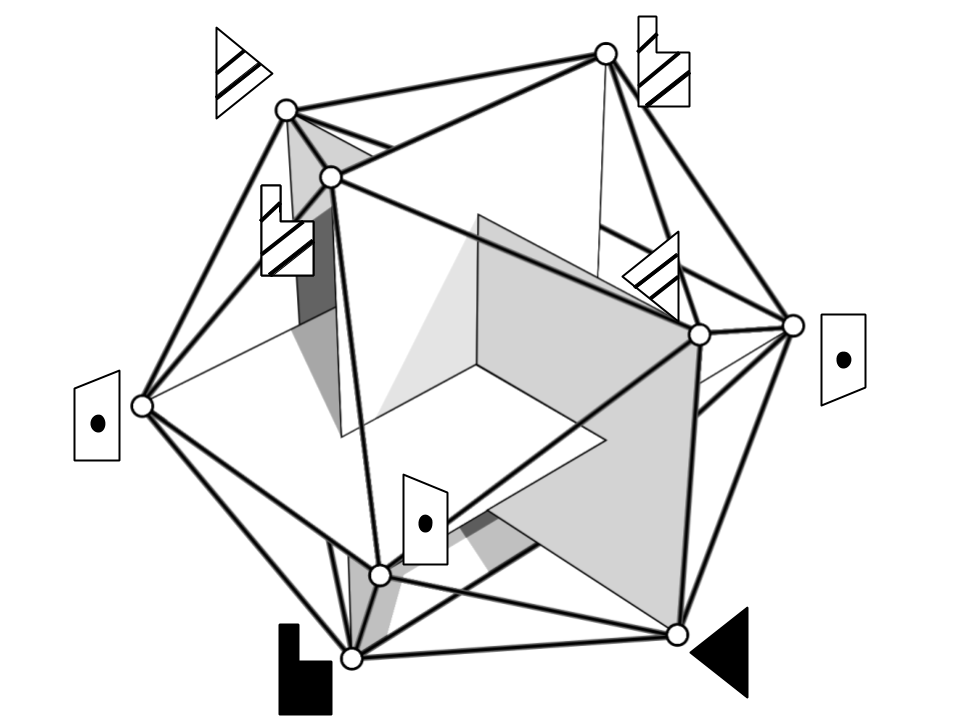}
  \caption{Directions on the Icosahedron}
\end{figure}

Using these spatial directions Laban established what he felt were natural sequences of movement inside the icosahedron and other spatial systems. He claimed these sequences followed natural patterns of the body and were governed by scientific and mathematical principles. Laban named many of these trace forms, including the defense and attack scales, and the girdle. We will start the exploration of the principles governing these scales on the most intriguing trace form, the primary scale on the icosahedron.  \\

\section*{The Primary Scale on the Icosahedron}

Laban established  what he called the standard (or primary) scale on the twelve vertices of the icosahedron, and it is displayed below. Laban recognized many symmetries and patterns used to establish this scale, including squares, triangles and diagonals in space. 
\begin{figure}[!h]
  \centering
  \includegraphics[scale=0.3]{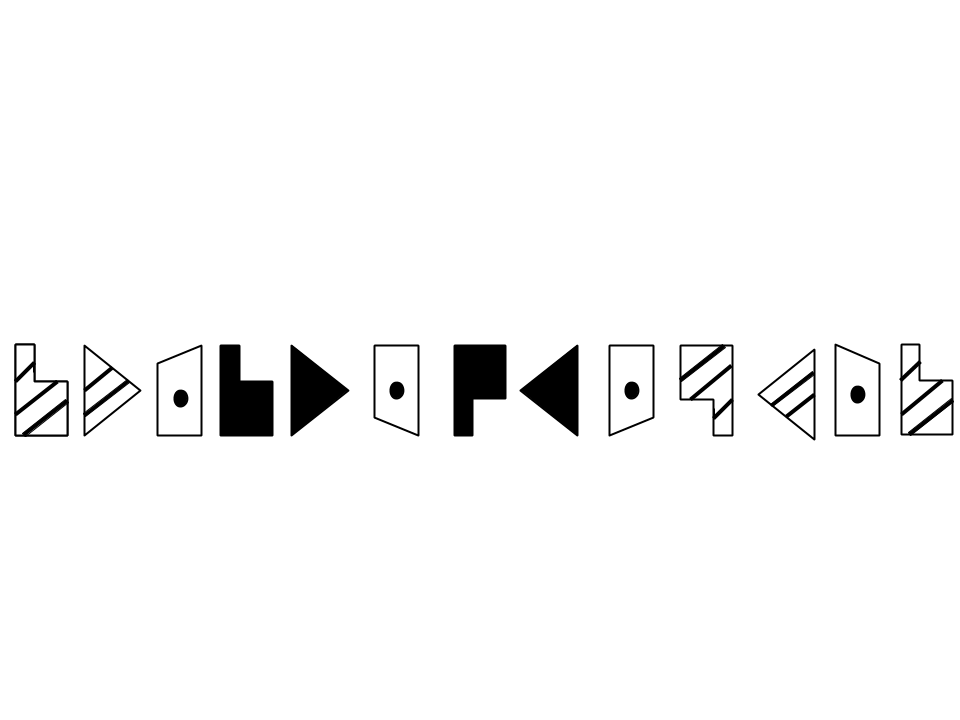}
 \caption{Laban's Primary Scale}
\end{figure}

\clearpage

Louis Horst was one of the first to recognize a possible connection between music composition and dance composition, and the devices of inversion and transposition in dance were translated from these same devices in music theory \cite{horst}. Taking a hint from music theory, in an attempt to better understand Laban's reasoning behind the primary scale, we can visualize the primary scale on a clock, akin to the musical clock from music theory \cite{mccartin}. \\

\begin{figure}[!h]
  \centering
  \includegraphics[scale=0.35]{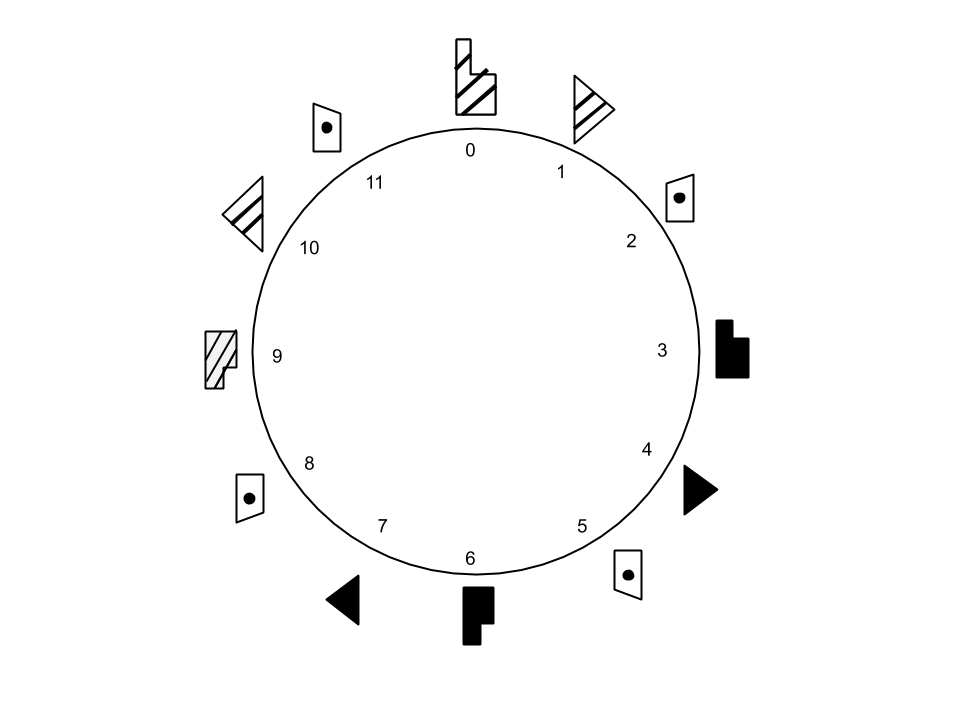}
  \caption{Primary Scale Clock}
\end{figure}

This illustrates how the primary scale can be represented by $Z_{12}$. And beautifully, the triangles, quadrangles and diameters referred to by Laban can be viewed as cosets of $Z_{12}$. The triangles correlate to $4Z_{12}$, the quadrangles $3Z_{12}$ and the diameters $6Z_{12}$. What Laban referred to as trace forms can be represented by paths connecting points of the clock. \\

\begin{figure}[!h]
  \centering
  \includegraphics[scale=0.35]{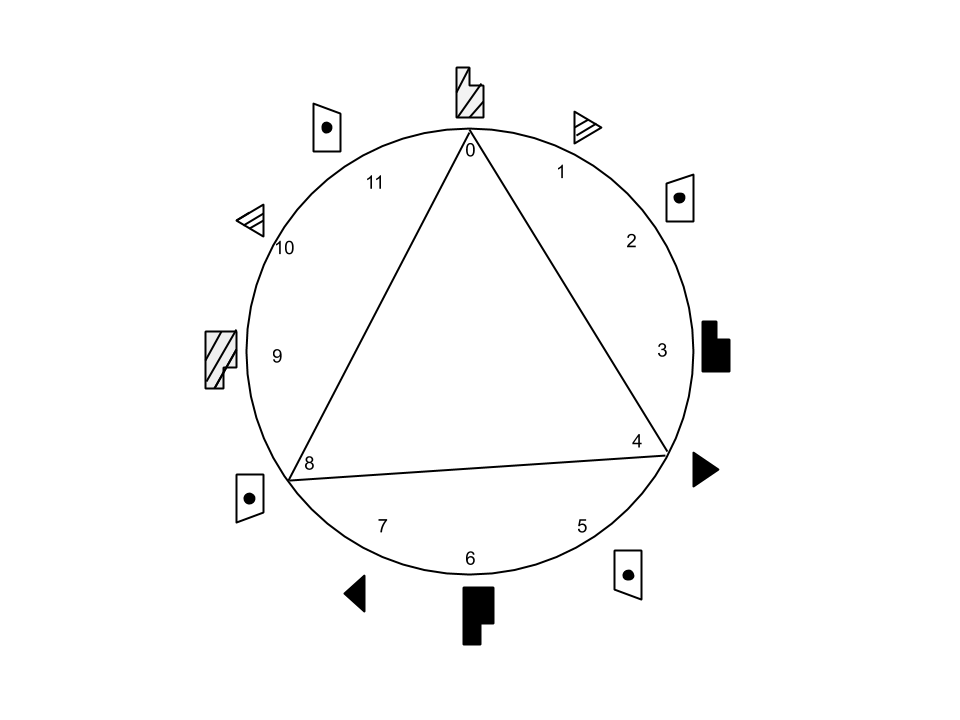}
  \caption{Triangles on the Clock}
\end{figure}

\begin{figure}[!h]
  \centering
  \includegraphics[scale=0.35]{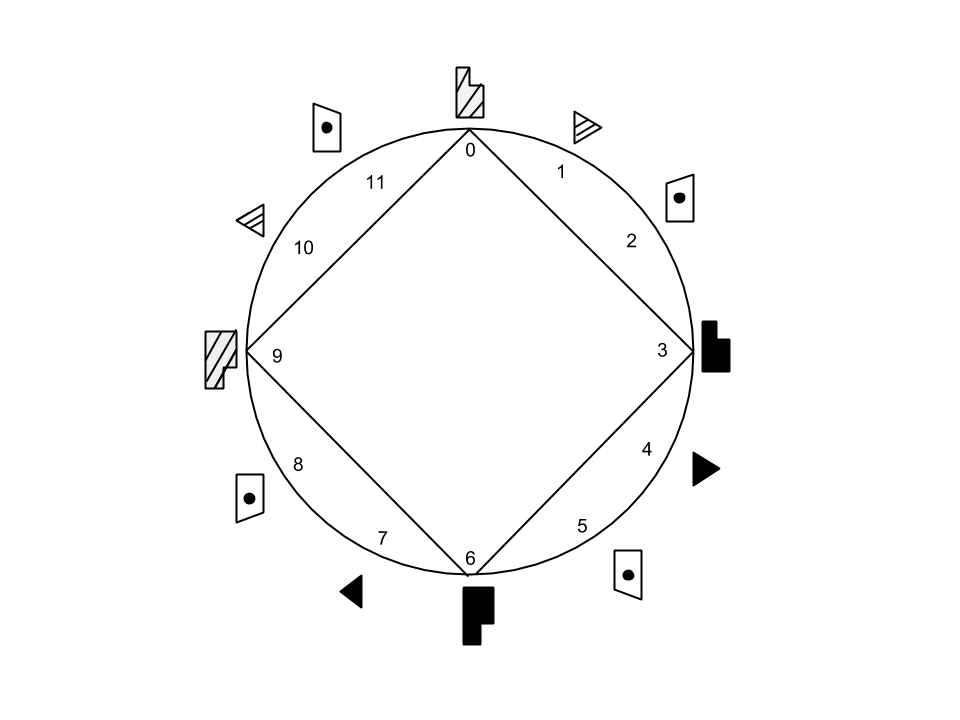}
  \caption{Squares on the Clock}
\end{figure}

\clearpage

\section*{Defining Inversion on the Primary Scale Clock}

Earlier we defined three different types of inversion in the icosahedron; front-back, left-right, and high low. We can now define another type of inversion, diametral inversion, on the primary scale clock. Diametral inversion could be defined as subsequent applications of the previously defined inversion. For instance the diametral inversion of an action in the right middle front direction would be an action in the left middle back direction, which could be obtained by applying a front back inversion, and then a left right inversion. This previous analysis focuses only on the final destination of the movement, rather than the path taken to arrive there, so we were only able to determine the inversion of a single movement, rather than the inversion of a sequence of movements. If we now refer only to diametral inversion, we can develop a method of determining the inversion for a whole sequence of movements using the primary scale clock.   

\begin{figure}[!h]
  \centering
  \includegraphics[scale=0.35]{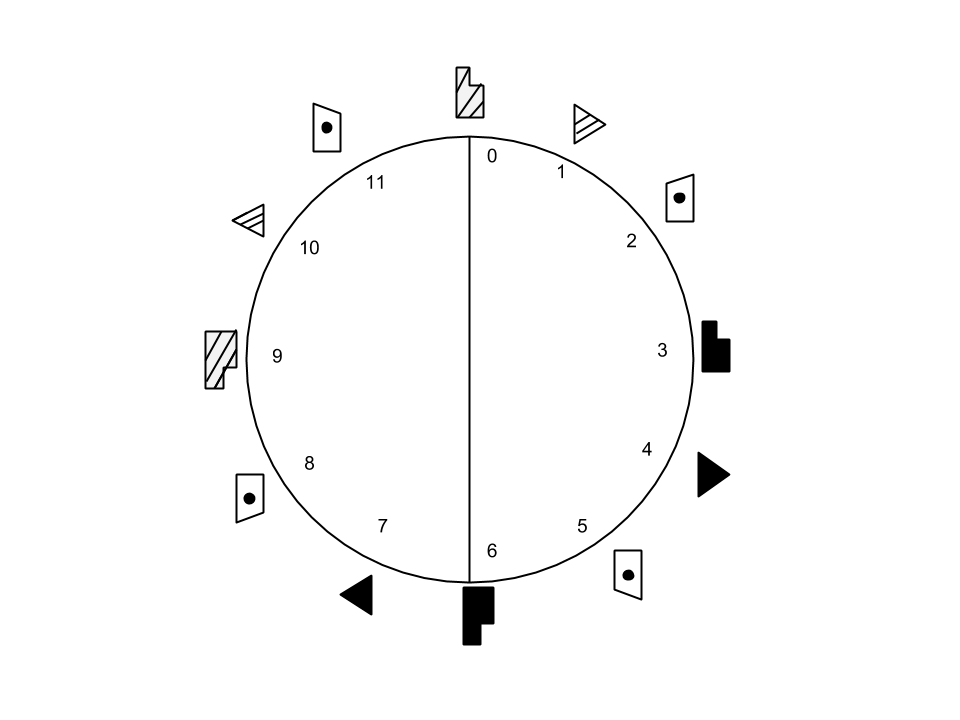}
  \caption{Diameters Defining Diametral Inversion}
\end{figure}

\begin{figure}[!h]
  \centering
  \includegraphics[scale=0.35]{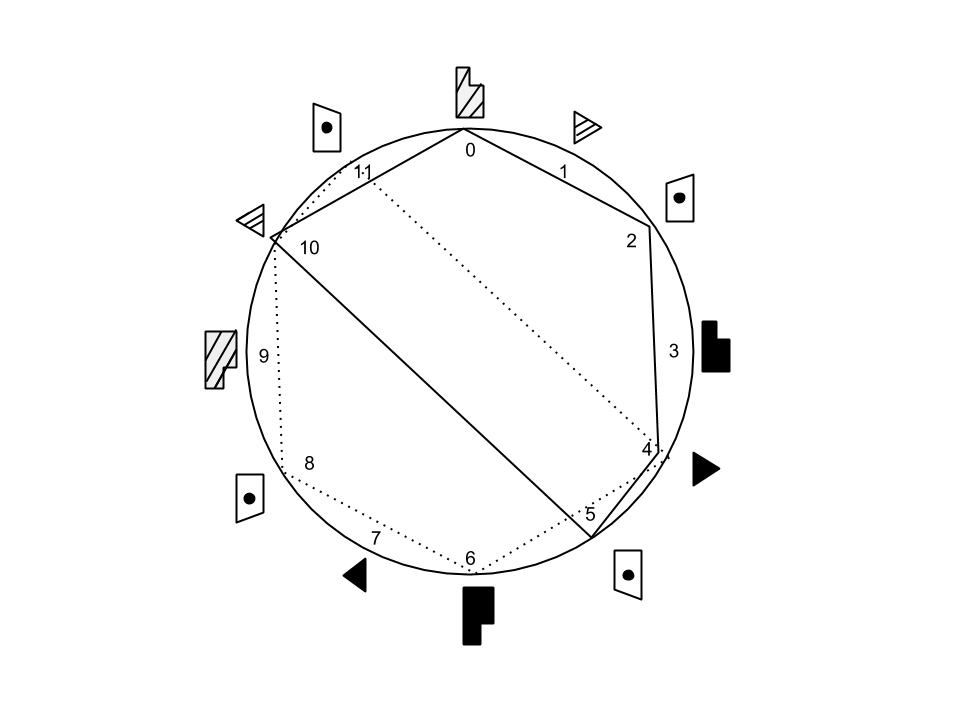}
  \caption{The "Girdle" Trace Form and It's Inversion}
\end{figure}

\clearpage
    
\section*{New Devices} The goal of this paper was to first mathematically establish the symmetries that define choreographic devices already in place, with the hope of then exploring additional symmetries to determine new devices. We have currently only looked at the concept of choreographic inversion, and shown how it can be determined by group actions on the octahedron and icosahedron. More analysis can be done on these polyhedra, and in the future to determine new devices it may be useful to examine non-regular and non-convex polyhedra as systems of spatial reference. \\

\end{document}